\documentclass[11pt]{article}
\usepackage{graphicx}
\usepackage{amsmath,amssymb,amsthm,amsfonts}
\usepackage[normal]{subfigure}
\usepackage{amssymb}
\usepackage{subfig}

\begin{document}
\title{Numerical Methods for Coupled Surface and Grain Boundary Motion}
\author{Zhenguo Pan , \  \  Brian Wetton \\  Department of Mathematics,
University of British Columbia,\\  Vancouver, B.C. Canada \\ }
\date{}
\maketitle
\begin{abstract}
We study the coupled surface and grain boundary motion in a
bicrystal in the context of the ``quarter loop" geometry. Two
types of physics motions are involved in this model: motion by
mean curvature and motion by surface diffusion. The goal is
finding a formulation that can describe the coupled motion and has
good numerical behavior when discretized. Two formulations are
proposed in this paper. One of them is given by a mixed order
parabolic system and the other is given by Partial Differential
Algebraic Equations. The parabolic formulation constitutes several
parabolic equations which model the two normal direction motions
separately. The performance of this formulation is good for a
short time simulation. It performs even better by adding an extra
term to adjust the tangential velocity of grid points. The PDAE
formulation preserves the scaled arc length property and performs
much better with no need to add an adjusting term. Both
formulations are proven to be well-posed in a simpler setting and
are solved by finite difference methods.
\end{abstract}

\vskip 0.2truein


\vskip 0.2truein
\section{Introduction}
Coupled surface and grain boundary motion is an important
phenomenon controlling the grain growth in materials processing
and synthesis. A commonly used model to study this coupled effect
is called ``quarter loop" geometry introduced by Dunn et
al.~\cite{dunn}.


\begin{figure}
\begin{center}
\includegraphics[width = 12cm, height = 5.1cm, clip]{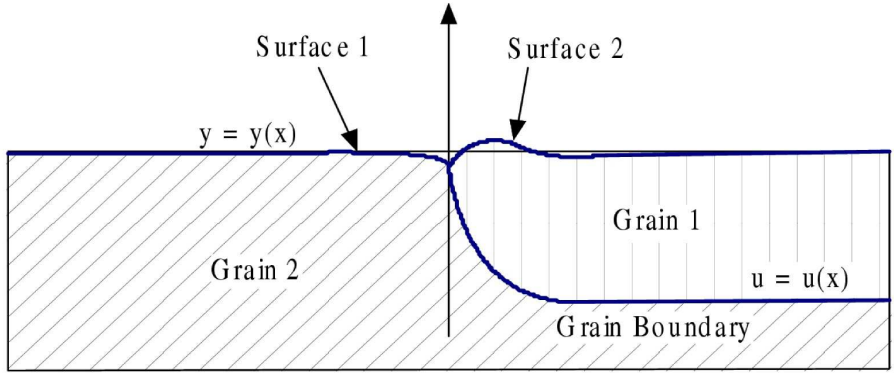}
\caption{The quarter loop bicrystal geometry.} \label{fig:fig1}
\end{center}
\end{figure}

In the quarter loop geometry, there are two grains between which
there is an interface called grain boundary as shown in
Fig.\ref{fig:fig1}. The two grains are of the same material and
differ only in their relative crystalline orientation. The grain
boundary runs parallel to a free surface before it turns up and
attaches to upper surfaces at a groove root. When heated at a
specific temperature, the grain boundary migrates to reduce the
surface energy and to heal the orientation mismatch. Since the
driving force is constant the grain boundary moves at a constant
velocity after a short time decay. It is reasonable to assume that
the bicrystal is uniform along the cross-section direction. Thus
it is reasonable for us to consider only two dimensional(2D)
geometry in this paper.

This geometry contains two types of motion. One of them is mean
curvature motion for the grain boundary. And the other one is
surface diffusion for the upper surfaces. More detail about this
model is give in ~\cite{amy}.

Motion by mean curvature is an evolution law in which the normal
velocity of an interface is proportional to its mean curvature.
More precisely, the motion of an interface $\Gamma$ satisfies
\begin{equation}
V_{c}=A\kappa
\end{equation}
Here $V_c$ denotes the velocity in the normal direction of
$\Gamma$, and $\kappa$ stands for the mean curvature of $\Gamma$.

First proposed by Mullins~\cite{mullins} to model the curvature
driven diffusion on the surface of a crystal, surface diffusion is
a different evolution law in which the normal velocity of an
interface is proportional to the surface Laplacian of mean
curvature. The motion of interface $\Gamma$ satisfies
\begin{equation}\label{eq:q1}
V_{d}=-B\Delta_s \kappa
\end{equation}
Here $V_{d}$ stands for the normal velocity and $\Delta_s$ stands
for the operator of surface Laplacian which is defined as
\begin{equation}
\Delta_s=\nabla_s\cdot\nabla_s \textrm{ where }
\nabla_s=\nabla-n\partial_n
\end{equation}

In two dimensions, surface diffusion can be reduced to a normal
direction motion with a speed function depending on the second
derivative of the curvature with respect to arc length, i.e.,
\begin{equation}\label{eq:surdiff}
V_{d}=-B\kappa_{ss}
\end{equation}
Here $s$ is arc length parametrization. Since the problem is
proposed in two dimensions, the motion by surface diffusion always
refers to equation (\ref{eq:surdiff}) instead of the general case
(\ref{eq:q1}) in this paper.

We shall prove later in the appendix that we could normalize $A$
and $B$ by rescaling the time and space. Since this reformulation
make the problem neither harder or easier we will take both $A$
and $B$ as one.

Surface diffusion is an intrinsically difficult problem to solve
numerically even in two dimensions. The main difficulty is that it
is stiff due to the fourth order derivatives and such that an
explicit time stepping strategy requires very small time steps.
Moreover, owing to the lack of a maximum principle, an embedded
curve may not stay embedded, in other words, it may become
self-intersected during the evolution.

Travelling wave solutions have been derived for the whole
nonlinear problem and for a linearized problem in~\cite{amy2}
and~\cite{amy} respectively. These analytic results are used to
verify numerical results in this paper.

\begin{figure}
\begin{center}
\includegraphics[width = 6cm, height = 4cm, clip]{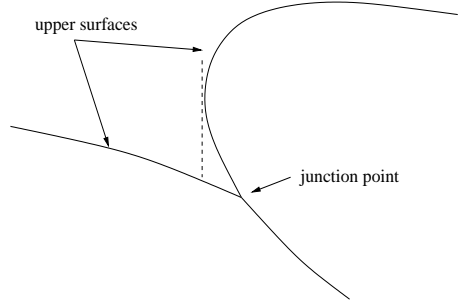}
\caption{An example when one of the surfaces is not
single-valued.} \label{fig:exception}
\end{center}
\end{figure}

The formulations in this paper will be proposed in parameterized
form. There are two reasons why we prefer the parameterized form.
Firstly, we try to set up a versatile formulation that is
extensible to other problems, such as the evolution of a closed
curve , other freely positioned triple junction problems and even
the evolution curve networks, wherever for a single type motion or
a mixed type motion. Secondly, even for the coupled grain boundary
motion, the function $y(x)$ which represents the exterior surface
may not be single-valued as shown in Fig.\ref{fig:exception} and
this phenomena is physically reasonable ~\cite{amy-1,amy0}. Also
for such consideration we will treat the exterior surface as two
curves separated by the triple junction in the following
discussion.

The outline of this paper is as follows. In section \ref{se:2},
parabolic equations are derived for the motion by mean curvature
and the motion by surface diffusion separately. The boundary
conditions including the triple junction conditions and domain
boundary conditions are discussed in section \ref{se:3} referring
to the analytic work of Novick-Cohen et al.~\cite{amy,amy2} and
Wong et al.~\cite{Donhong}. In section \ref{se:4}, the
well-posedness for a linear parabolic system that is closely
related to the full nonlinear problem is analyzed and followed by
a discussion about the artificial tangential condition. From
section \ref{se:5} to \ref{se:tangential} we discuss the numerical
details including discretization, time stepping and some other
numerical issues. We start the discussion for a PDAE system in
section \ref{se:8}. An interesting computational example of
surface diffusion is given in section \ref{se:star}. The linear
well-posedness of the PDAE system is analyzed in section
\ref{se:9}.

\section{Cartesian Formulation}
In this section, we consider the problem in the cartesian
coordinate system which is give as below(see~\cite{amy}).
\begin{eqnarray}\label{eq:car}
y_t&=&-\Big[\frac{1}{(1+y_x^2)^{1/2}}\Big[\frac{y_{xx}}{(1+y_x^2)^{3/2}}\Big]_x\Big]_x,\quad
t>0,x\in(-\infty,s(t)^-)\cup(s(t)^+,\infty),\nonumber\\
u_t&=&u_{xx}(1+u_x^2)^{-1},\quad t>0,x>s(t)
\end{eqnarray}
with following triple junction conditions
\begin{eqnarray}\label{eq:carcon}
&&y(s(t)^+,t)=y(s(t)^-,t)=u(s(t)^+,t),\quad t>0\nonumber\\
&&\arctan(y_x(s(t)^+,t))-\arctan(y_x(s(t)^-,t))=2\arcsin(\frac{\gamma_{grain
}}{2\gamma_{exterior}})\nonumber\\
&&\arctan(u_x(s(t)^+,t))=-\frac{\pi}{2}+\frac{1}{2}[\arctan(y_x(s(t)^+,t))+\arctan(y_x(s(t)^-,t))],\,t>0\\
&&\frac{y_{xx}}{(1+y_x^2)^{3/2}}\Big|_{(s(t)^+,t)}=\frac{y_{xx}}{(1+y_x^2)^{3/2}}\Big|_{(s(t)^-,t)},\,t>0\nonumber\\
&&\Big[\frac{1}{(1+y_x^2)^{1/2}}\Big[\frac{y_{xx}}{(1+y_x^2)^{3/2}}\Big]_x\Big]_{(s(t)^+,t)}=\Big[\frac{1}{(1+y_x^2)^{1/2}}\Big[\frac{y_{xx}}{(1+y_x^2)^{3/2}}\Big]_x\Big]_{(s(t)^-,t)},\,t>0\nonumber\\
&&y(+\infty,t)=y(-\infty,t)=0,\, t>0\nonumber\\
&&u(+\infty,t)=-1,\,t>0\nonumber
\end{eqnarray}
Here $y=y(x,t)$ stands for the height of the two exterior
surfaces. $u=u(x,t)$ is the height of the grain boundary and
$s(t)$ denotes the location of the junction where the three
surfaces meet.

Since the junction is moving it is not straightforward to solve
this system numerically. We fix the junction by making the
following transform,
\begin{eqnarray}
\bar{x}=x-s(t)
\end{eqnarray}
We let $y(\bar{x},t)=y(x,t)$ and $u(\bar{x},t)=u(x,t)$. Therefor,
\begin{eqnarray*}
y_x(x,t)&=&y_{\bar{x}}(\bar{x},t)\\
y_t(x,t)&=&y_t(\bar{x},t)-y_{\bar{x}}(\bar{x},t)s_t\\
u_x(x,t)&=&u_{\bar{x}}(\bar{x},t)\\
u_t(x,t)&=&u_t(\bar{x},t)-u_{\bar{x}}(\bar{x},t)s_t
\end{eqnarray*}
Then system (\ref{eq:car}) becomes
\begin{eqnarray}\label{eq:car2}
y_t&=&-\Big[\frac{1}{(1+y^2_{\bar{x}})^{1/2}}\Big[\frac{y_{\bar{x}\bar{x}}}{(1+y_{\bar{x}}^2)^{3/2}}\Big]_{\bar{x}}\Big]_{\bar{x}}+y_{\bar{x}}(\bar{x},t)s_t,\quad
t>0,\bar{x}\in(-\infty,0)\cup(0,\infty),\nonumber\\
u_t&=&u_{\bar{x}\bar{x}}(1+u_{\bar{x}}^2)^{-1}+y_{\bar{x}}(\bar{x},t)s_t,\quad
t>0,\bar{x}>0
\end{eqnarray}
And the boundary conditions (\ref{eq:carcon}) become
\begin{eqnarray}\label{eq:carcon2}
&&y(0^+,t)=y(0^-,t)=u(0^+,t),\quad t>0\nonumber\\
&&\arctan(y_{\bar{x}}(0^+,t))-\arctan(y_{\bar{x}}(0^-,t))=2\arcsin(\frac{\gamma_{grain
}}{2\gamma_{exterior}})\nonumber\\
&&\arctan(u_{\bar{x}}(0^+,t))=-\frac{\pi}{2}+\frac{1}{2}[\arctan(y_{\bar{x}}(0^+,t))+\arctan(y_{\bar{x}}(0^-,t))],\,t>0\\
&&\frac{y_{{\bar{x}}{\bar{x}}}}{(1+y_{\bar{x}}^2)^{3/2}}\Big|_{(0^+,t)}=\frac{y_{{\bar{x}}{\bar{x}}}}{(1+y_{\bar{x}}^2)^{3/2}}\Big|_{(0^-,t)},\,t>0\nonumber\\
&&\Big[\frac{1}{(1+y_{\bar{x}}^2)^{1/2}}\Big[\frac{y_{{\bar{x}}{\bar{x}}}}{(1+y_{\bar{x}}^2)^{3/2}}\Big]_{\bar{x}}\Big]_{(0^+,t)}=\Big[\frac{1}{(1+y_{\bar{x}}^2)^{1/2}}\Big[\frac{y_{{\bar{x}}{\bar{x}}}}{(1+y_{\bar{x}}^2)^{3/2}}\Big]_{\bar{x}}\Big]_{(0^-,t)},\,t>0\nonumber\\
&&y(+\infty,t)=y(-\infty,t)=0,\, t>0\nonumber\\
&&u(+\infty,t)=-1,\,t>0\nonumber
\end{eqnarray}

\begin{figure}
\begin{center}
\includegraphics[width = 7cm, height = 6cm, clip]{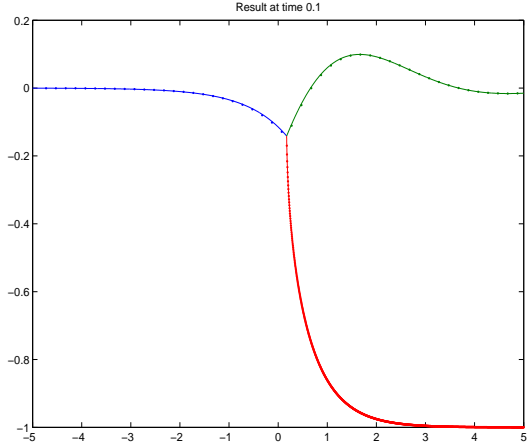}
\caption{Numerical result for system
(\ref{eq:car2})-(\ref{eq:carcon2}) with $m=0.5$. Dotted line:
numerical result; Solid line: travelling wave solution.}
\label{fig:car}
\end{center}
\end{figure}

This system could easily be discretized using standard finite
difference schemes on a fixed staggered grid. A numerical result
is shown in Fig.\ref{fig:car}.

The disadvantage of using cartesian formulation is that it is not
applicable to non-single valued case as shown in
Fig.\ref{fig:exception}. What's more, since the grain boundary is
nearly singular at the junction, it requires very small grid size
for accuracy. For wider application we consider two parametric
formulations in the rest of this paper.

\section{A Parabolic Formulation}\label{se:2}
In this section we derive a parabolic system to model the coupled
motion. Here and throughout this paper we use
$X=(u(\cdot),v(\cdot))$ to represent a parameterized curve with
$u(\cdot)$ and $v(\cdot)$ being the coordinates.

Several more notations should be introduced as well. The arc
length parameter is denoted by $X(s)=X(u(s),v(s))$ and any other
parameter is denoted by $X(\sigma)=X(u(\sigma),v(\sigma))$.
$\vec{t}$ and $\vec{n}$ stand for unit tangential direction and
unit normal direction respectively. $\kappa$ stands for curvature.
Although all final equations are parameterized by $\sigma$, the
arc length parametrization is useful for the intermediate
deviations.
\subsection{Motion by Mean Curvature}
We first derive a parabolic equation to describe the motion by
mean curvature. Similar discussion has been addressed
in~\cite{brian} and~\cite{garcke}. We give a brief description for
reader's convenience.

With the notations introduced above one has
\begin{eqnarray*}
X_s=\vec{t}\nonumber
\end{eqnarray*}
\begin{eqnarray*}
X_{ss}=\kappa\vec{n}
\end{eqnarray*}
Here the subscript $s$ stands for the derivative of $X$ with
respect to arc length $s$. Direct computation shows that
\begin{equation}\label{eq:8}
X_{\sigma}=X_s\frac{dS(\sigma)}{d\sigma}=X_s|X_\sigma|
\end{equation}
where $S(\sigma)$ is defined by
\begin{eqnarray}\label{eq:9}
S(\sigma)=\int_{\sigma_0}^{\sigma}\sqrt{u_{\sigma}^2+v_{\sigma}^2}d\sigma
\end{eqnarray}
which stands for the length of the curve from point $X(\sigma_0)$
to $X(\sigma)$. $|X_\sigma|$ is $L_2$ norm of $X_\sigma$ defined
by
\begin{eqnarray*}
|X_\sigma|=\sqrt{u_\sigma^2+v_\sigma^2}
\end{eqnarray*}
 Differentiate equation
(\ref{eq:8}) with respect to $\sigma$ to botain
\begin{equation}
X_{\sigma\sigma}=X_{ss}|X_\sigma|^2+X_s|X_\sigma|_s|X_\sigma|
\end{equation}
By previous derivations one can compute the normal component of
vector $\frac{X_{\sigma\sigma}}{|X_\sigma|^2}$ and obtains
\begin{eqnarray}\label{eq:normal}
\frac{X_{\sigma\sigma}}{|X_\sigma|^2}\cdot\vec{n}&=&X_{ss}\cdot
\frac{X_{ss}}{\kappa}+\frac{|X_\sigma|_s}{|X_\sigma|}X_s\cdot\frac{X_{ss}}{\kappa}\nonumber\\
&=& \kappa
\end{eqnarray}
Thus, if we set up a formulation:
\begin{equation}\label{eq:curvature1}
X_t=\frac{X_{\sigma\sigma}}{|X_\sigma|^2}
\end{equation}
it is obvious that the motion described by (\ref{eq:curvature1})
has normal velocity $\kappa$. This gives us an option to describe
the motion by mean curvature. Equation (\ref{eq:curvature1}) is
fully parabolic which means it is parabolic in both the normal
component and tangential component.

There are some other equations that can also describe curvature
motion, for example,
\begin{eqnarray}\label{eq:alternate}
X_t=\kappa \vec{n}
\end{eqnarray}
But this system is not fully parabolic. A linearization shows that
this system is parabolic in the normal component and hyperbolic in
the tangential component. A discretization of system
(\ref{eq:alternate}) will not have the good numerical properties
as those of a fully parabolic system due to the lack of regularity
in the parametrization as shown in~\cite{brian}.

\subsection{Motion by Surface Diffusion} Considering the good
properties of a parabolic formulation, we hope to find a parabolic
formulation for motion by surface diffusion. By analogy with the
approach to the mean curvature motion described above, we try the
following form:
\begin{equation}\label{eq:14}
X_t=-\frac{X_{\sigma\sigma\sigma\sigma}}{|X_\sigma|^4}+
L(X_{\sigma\sigma\sigma},X_{\sigma\sigma},X_{\sigma})
\end{equation}
where $L(X_{\sigma\sigma\sigma},X_{\sigma\sigma},X_{\sigma})$
includes some lower order terms and will be determined such that
\begin{eqnarray}\label{eq:145}
X_t\cdot \vec
{n}=(-\frac{X_{\sigma\sigma\sigma\sigma}}{|X_\sigma|^4}+
L(X_{\sigma\sigma\sigma},X_{\sigma\sigma},X_{\sigma}))\cdot
\vec{n}=-\kappa_{ss}
\end{eqnarray}

We focus our study on finding out
$L(X_{\sigma\sigma\sigma},X_{\sigma\sigma},X_{\sigma})$ in the
rest of this section.

Note first the following equation (see appendix for proof),
\begin{equation}\label{eq:ksss}
(X_{ssss}+\kappa^2X_{ss})\cdot \vec{n}=\kappa_{ss}
\end{equation}
Compare equation (\ref{eq:145}) and (\ref{eq:ksss}) to get a
choice for $L$,
\begin{equation}\label{eq:newl}
L(X_{\sigma\sigma\sigma},X_{\sigma\sigma},X_{\sigma})=\frac{X_{\sigma\sigma\sigma\sigma}}{|X_\sigma|^4}-X_{ssss}
-\kappa^2 X_{ss}
\end{equation}
One will find later that the fourth order terms appeared in
(\ref{eq:newl}) could be cancelled with each other and such that
$L$ involves only third or lower order derivatives.

Start by equation (\ref{eq:8}) and differentiate several times
with respect to $\sigma$ to get following relations,
\begin{equation}\label{eq:17}
X_{\sigma\sigma}=X_{ss}S_\sigma^2+X_sS_{\sigma\sigma}
\end{equation}
\begin{equation}\label{eq:15}
X_{\sigma\sigma\sigma}=X_{sss}S_\sigma^3+3X_{ss}S_\sigma
S_{\sigma\sigma}+X_sS_{\sigma\sigma\sigma}
\end{equation}
\begin{equation}\label{eq:16}
X_{\sigma\sigma\sigma\sigma}=X_{ssss}S_\sigma^4+6X_{sss}S_\sigma^2S_{\sigma\sigma}+4X_{ss}S_\sigma
S_{\sigma\sigma\sigma}+3X_{ss}S_{\sigma\sigma}^2+X_sS_{\sigma\sigma\sigma\sigma}
\end{equation}
Dividing through equation (\ref{eq:16}) by $S_\sigma^4$ and
noticing the fact that $|X_\sigma|=S_\sigma$ one obtains
\begin{equation}\label{eq:116}
\frac{X_{\sigma\sigma\sigma\sigma}}{|X_\sigma|^4}=X_{ssss}+6\frac{S_{\sigma\sigma}}{S_\sigma^2}X_{sss}
+4\frac{S_{\sigma\sigma\sigma}}{S_\sigma^3}X_{ss}+3\frac{S_{\sigma\sigma}^2}{S_\sigma^4}X_{ss}
+\frac{S_{\sigma\sigma\sigma\sigma}}{S_\sigma^4}X_s
\end{equation}
Substitute equation (\ref{eq:116}) into (\ref{eq:newl}), rewrite
arc length parametrization $s$ into $\sigma$ using
(\ref{eq:17})-(\ref{eq:15}). Since vector $X_s$ is perpendicular
to $\vec{n}$ and has no contribution to normal direction we can
ignore all $X_s$ terms and obtain
\begin{equation}\label{eq:F}
L(X_{\sigma\sigma\sigma},X_{\sigma\sigma},X_{\sigma})=
6\frac{S_{\sigma\sigma}}{S_\sigma^2}\frac{X_{\sigma\sigma\sigma}}{|X_\sigma|^3}-
15\frac{S_{\sigma\sigma}^2}{S_\sigma^4}\frac{X_{\sigma\sigma}}{|X_\sigma|^2}+
4\frac{S_{\sigma\sigma\sigma}}{S_\sigma^3}\frac{X_{\sigma\sigma}}{|X_\sigma|^2}
-\kappa^2\frac{X_{\sigma\sigma}}{|X_\sigma|^2}
\end{equation}

Substitute equation (\ref{eq:F}) back into (\ref{eq:14}) and
collect to get the scheme as
\begin{equation}\label{eq:scheme}
X_t=-\frac{X_{\sigma\sigma\sigma\sigma}}{|X_\sigma|^4}+
6S_{\sigma\sigma}\frac{X_{\sigma\sigma\sigma}}{|X_\sigma|^5}-
(15\frac{S_{\sigma\sigma}^2}{|X_\sigma|^4}-
4\frac{S_{\sigma\sigma\sigma}}{|X_\sigma|^3}+\kappa^2)\frac{X_{\sigma\sigma}}{|X_\sigma|^2}
\end{equation}

We would like to point out that the choice of $L$ is not unique. A
similar expression has been given by Garcke et al.
in~\cite{garcke}.

\subsection{The Parabolic System}
We now give the fully parabolic system,
 \begin{eqnarray}
X_{t}^1&=&\frac{X^1_{\sigma\sigma}}{|X^1_{\sigma}|^2}\nonumber\\
X_{t}^2&=&-\frac{X^2_{\sigma\sigma\sigma\sigma}}{|X^2_{\sigma}|^4}+
6S_{\sigma\sigma}\frac{X^2_{\sigma\sigma\sigma}}{|X^2_{\sigma}|^5}-
(15\frac{S_{\sigma\sigma}^2}{|X^2_{\sigma}|^4}-
4\frac{S_{\sigma\sigma\sigma}}{|X^2_{\sigma}|^3}+\kappa^2)\frac{X^2_{\sigma\sigma}}{|X^2_{\sigma}|^2}\label{eq:sys}\\
X^3_{t}&=&-\frac{X^3_{\sigma\sigma\sigma\sigma}}{|X^3_{\sigma}|^4}+
6S_{\sigma\sigma}\frac{X^3_{\sigma\sigma\sigma}}{|X^3_{\sigma}|^5}-
(15\frac{S_{\sigma\sigma}^2}{|X^3_{\sigma}|^4}-
4\frac{S_{\sigma\sigma\sigma}}{|X^3_{\sigma}|^3}+\kappa^2)\frac{X^3_{\sigma\sigma}}{|X^3_{\sigma}|^2}\nonumber
\end{eqnarray}
where $X^1$ stands for the grain boundary and $X^2,X^3$ stand for
the left branch and right branch of the upper surface
respectively. All curves are represented by $X(\sigma)$ with
$\sigma\in[0,\infty)$.

This system will be solved numerically with the boundary
conditions discussed in the next section.

\section{Boundary Conditions}\label{se:3}
The grain boundary and the two upper surfaces meet together at one
end which is referred as triple junction. The other end of the
three curves tends to infinity in the quarter loop geometry. For
numerical reasons, we compute this problem in a bounded domain.
This domain is chosen large enough such that it can simulate the
motion at least for a short time. This restriction is reasonable
since the curves are asymptotically flat for the parts far away
from the triple junction. All computations presented in this paper
are constrained in a finite domain $[-6,12]$ and the curves are
parameterized with $\sigma\in [0,1]$.

At $\sigma=0$ the three curves meet at a triple junction and at
$\sigma=1$ the three curves meet the artificial domain boundary
separately.

\subsection{Triple Junction Conditions at $\sigma=0$}\label{se:junction}
\begin{figure}
\begin{center}
\includegraphics[width = 8.5cm, height = 5cm, clip]{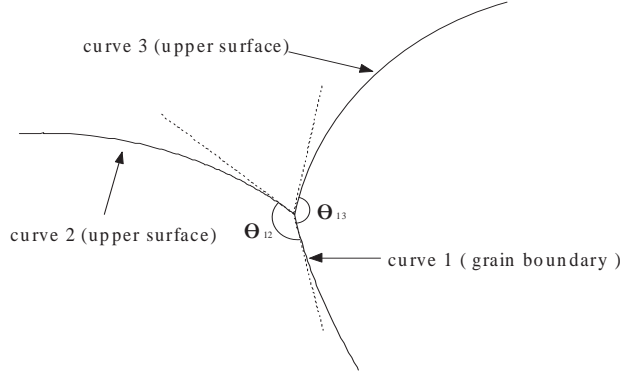}
\caption{Sketch of the grain boundary groove.} \label{fig:fig2}
\end{center}
\end{figure}
We first discuss the boundary conditions at the triple junction.
First of all, three curves should have common coordinates at
$\sigma=0$, i.e.,
\begin{equation}\label{eq:cond1}
X^1(0,t)=X^2(0,t)=X^3(0,t)
\end{equation}

 By Young's law we have two more conditions which are referred as angle
 conditions:
\begin{eqnarray}\label{eq:cond2}
&&\frac{X^1_{\sigma}}{|X^1_{\sigma}|}\cdot \frac{X^2_{\sigma}}{|X^2_{\sigma}|}=\cos\theta_{12}=\cos(\frac{\pi}{2}+\arcsin\frac{m}{2})\label{eq:cond3}\\
&&\frac{X^1_{\sigma}}{|X^1_{\sigma}|}\cdot
\frac{X^3_{\sigma}}{|X^3_{\sigma}|}=\cos\theta_{13}=\cos(\frac{\pi}{2}+\arcsin\frac{m}{2})\label{eq:cond3}
\end{eqnarray}
where $\theta_{ij}$ denotes the angle between curve $i$, $j$ and
$m=\gamma_{grain }/{\gamma_{exterior}}$ is a constant measuring
the relative surface energy between the grain boundary and
exterior surface.

The continuity of the surface chemical potentials implies that
\begin{equation}\label{eq:potentials}
\kappa^2=-\kappa^3\qquad
(\kappa=\frac{X_{\sigma\sigma}}{|X_\sigma|^2}\cdot
\frac{X_\sigma^\bot}{|X_\sigma|})
\end{equation}
Here the superscripts are indices of curves.

And the balance of mass flux implies that
\begin{equation}\label{eq:mass}
\kappa^2_{s}=\kappa^3_{s}\qquad
(\kappa_s=\frac{X_{\sigma\sigma\sigma}\cdot
X_\sigma^\perp}{|X_\sigma|^4}-
3\frac{|X_\sigma|_\sigma(X_{\sigma\sigma}\cdot X_\sigma^\perp)}
{|X_\sigma|^5})
\end{equation}
where the expression for $\kappa_s$ is obtained by taking the
derivative of the expression of $\kappa$ directly.

We must be careful about condition (\ref{eq:potentials}).
Basically, we need the two upper surfaces have the same convexity.
Since $\sigma$ has opposite directions for the two curves the odd
time derivatives will have opposite signs when computed by
parametric form. Thus we should put a minus sign for
(\ref{eq:potentials}) and keep the same for (\ref{eq:mass}).

\subsection{Boundary Conditions at $\sigma=1$}
At the other ends of the curves we put several artificial
conditions such that they do not move during evolution and keep
being flat. This is reasonable since they start being flat and
they will not be influenced by the motion of the triple junction
in a short time. The following conditions are imposed at
$\sigma=1$,
\begin{eqnarray*}
&& X_{t}^i(1,t)=\textbf{0}\qquad \textrm{for}\, i=1,2,3\\
&& X_{\sigma\sigma}^i(1,t)=\textbf{0}\qquad \textrm{for}\, i=2,3
\end{eqnarray*}
\subsection{Artificial Tangential Conditions}
We should point out that the whole system contains two second
order equations and four fourth order equations and it should have
ten conditions at the junction point for well-posedness. Recall
that there are only eight junction conditions as have been
addressed above. Thus, we need two more conditions. There are
several options to impose the extra conditions. And we will prove
later that different conditions could only change the
parametrization of the curves and will not change the profiles of
the curves. Since these conditions do change the tangential
velocities of the grid nodes we refer them as artificial
tangential conditions. As one of the options, the following two
conditions are applied into the system:
\begin{equation}
X_{\sigma\sigma}^i\cdot X_{\sigma}^i=0 \,\,\quad \textrm{for}
\,i=2,3
\end{equation}
\section{Well-posedness for the Parabolic System}\label{se:4}
In this section, we analyze the well-posedness of the system
proposed above. We linearize around fixed straight line solutions
and get a system that has the same highest order parabolic
behavior as the original problem. The well-posedness we do gives
the conditions that match those that in more complicated nonlinear
analysis gives, where such analysis exists. And therefore, we
believe the results of the analysis should apply to the full
nonlinear problem.
\subsection{Linearization of the System}\label{se:linearized}
To linearize the system we consider a perturbation expansion
around the tangential direction at the triple junction for each
curve, i.e.,
\begin{eqnarray*}
&&X^1=d_1\sigma+\epsilon \bar{X}^1+O(\epsilon^2)\\
&&X^2=d_2\sigma+\epsilon \bar{X}^2+O(\epsilon^2)\\
&&X^3=d_3\sigma+\epsilon \bar{X}^3+O(\epsilon^2)
\end{eqnarray*}
where $d_i=(d_{i1},d_{i2})$ is a constant vector standing for the
unit tangential direction. Substitute above equations into
(\ref{eq:sys}), linearize and keep the leading order terms to get
a linear system:
\begin{eqnarray}
&&\bar{X}^1_{t}=\bar{X}^1_{\sigma\sigma}\nonumber\\
&&\bar{X}^2_{t}=-\bar{X}^2_{\sigma\sigma\sigma\sigma}\label{eq:linearized}\\
&&\bar{X}^3_{t}=-\bar{X}^3_{\sigma\sigma\sigma\sigma}\nonumber
\end{eqnarray}
For convenience, we omit the bar above $X$ in following
discussion.

The linearization of the triple junction conditions is
straightforward.
\begin{itemize}
\item Common point at $\sigma=0$:
\begin{eqnarray*}
X^1=X^2=X^3
\end{eqnarray*}
\item Angle conditions:
\begin{eqnarray*}
d_1\cdot X^2_{\sigma}+d_2\cdot X^1_{\sigma}-(d_1\cdot
d_2)(d_1\cdot
X^1_{\sigma}+d_2\cdot X^2_{\sigma})=0\\
d_1\cdot X^3_{\sigma}+d_3\cdot X^1_{\sigma}-(d_1\cdot
d_3)(d_1\cdot X^1_{\sigma}+d_3\cdot X^3_{\sigma})=0
\end{eqnarray*}
\item Continuity of surface chemical potentials:
\begin{eqnarray*}
X^2_{\sigma\sigma}\cdot d_2^\perp =-X^3_{\sigma\sigma} \cdot
d_3^\perp
\end{eqnarray*}
\item Balance of mass flux:
\begin{eqnarray*}
X^2_{\sigma\sigma\sigma}\cdot d_2^\perp =X^3_{\sigma\sigma\sigma}
\cdot d_3^\perp
\end{eqnarray*}
\item Artificial tangential conditions:
\begin{eqnarray*}
X^2_{\sigma\sigma}\cdot d_2=0\\
X^3_{\sigma\sigma}\cdot d_3=0
\end{eqnarray*}
\end{itemize}
The linear system (\ref{eq:linearized}) can be solved using
Laplace transforms to get {\large
\begin{eqnarray}\label{eq:sol}
\left\{\begin{array}{l}
                    u_1=A_{11}e^{-\sqrt{s}\sigma}\\
                    v_1=A_{12}e^{-\sqrt{s}\sigma}\\
                    u_2=A_{21}e^{\lambda_1\sigma}+B_{21}e^{\lambda_2\sigma}\\
                    v_2=A_{22}e^{\lambda_1\sigma}+B_{22}e^{\lambda_2\sigma}\\
                    u_3=A_{31}e^{\lambda_1\sigma}+B_{31}e^{\lambda_2\sigma}\\
                    v_3=A_{32}e^{\lambda_1\sigma}+B_{32}e^{\lambda_2\sigma}\\
                    \end{array}
                    \right.
\end{eqnarray}}
where
\begin{eqnarray*}
\lambda_1=(-\frac{\sqrt{2}}{2}+\frac{\sqrt{2}}{2}i)\sqrt[4]{s}\quad
\lambda_2=(-\frac{\sqrt{2}}{2}-\frac{\sqrt{2}}{2}i)\sqrt[4]{s}
\end{eqnarray*}
and here $s$ temporally stands for the transformed time variable
of Laplace transform.

For simplicity, we first suppose the angles between any two curves
are $\frac{2}{3}\pi$. Substituting solution (\ref{eq:sol}) into
boundary conditions one obtains a $10 \times 10$ coefficient
matrix $M$(transposed)
\begin{eqnarray}\label{eq:matrix}
\left(\begin{array}{cccccccccc} 1&0&0&0&(-d_{21}-\frac{1}{2}
d_{11})\sqrt{s}&(-d_{31}-\frac{1}{2}
d_{11})\sqrt{s}&0&0&0&0\\
0&0&1&0&(-d_{22}-\frac{1}{2}d_{12})\sqrt{s}&(-d_{32}-\frac{1}{2}d_{12})\sqrt{s}&0&0&0&0\\
-1&1&0&0&(d_{11}+\frac{1}{2}
d_{21})\lambda_1&0&-d_{22}\lambda_1^2&-d_{22}\lambda_1^3&d_{21}\lambda_1^2&0\\
-1&1&0&0&(d_{11}+\frac{1}{2}
d_{21})\lambda_2&0&-d_{22}\lambda_2^2&-d_{22}\lambda_2^3&d_{21}\lambda_2^2&0\\
0&0&-1&1&(d_{12}+\frac{1}{2}
d_{22})\lambda_1&0&d_{21}\lambda_1^2&d_{21}\lambda_1^3&d_{22}\lambda_1^2&0\\
0&0&-1&1&(d_{12}+\frac{1}{2}d_{22})\lambda_2&0&d_{21}\lambda_2^2&d_{21}\lambda_2^3&d_{22}\lambda_2^2&0\\
 0&-1&0&0&0&(d_{11}+\frac{1}{2}
 d_{31})\lambda_1&-d_{32}\lambda_1^2&d_{32}\lambda_1^3&0&d_{31}\lambda_1^2\\
0&-1&0&0&0&(d_{11}+\frac{1}{2}
d_{31})\lambda_2&-d_{32}\lambda_2^2&d_{32}\lambda_2^3&0&d_{31}\lambda_2^2\\
0&0&0&-1&0&(d_{12}+\frac{1}{2}
d_{32})\lambda_1&d_{31}\lambda_1^2&-d_{31}\lambda_1^3&0&d_{32}\lambda_1^2\\
0&0&0&-1&0&(d_{12}+\frac{1}{2}d_{32})\lambda_2&d_{31}\lambda_2^2&-d_{31}\lambda_2^3&0&d_{32}\lambda_2^2\\
 \end{array}
\right)
\end{eqnarray}

Linear well-posedness requires that the determinant of matrix $M$
is nonsingular for any $s$ satisfying Re$(s)>0$. Since the
well-posedness depends only on their relative positions, we
suppose further that
\begin{eqnarray*}
d_1=\left(\begin{array}{c} 0\\\\-1\end{array}\right)\quad
d_2=\left(\begin{array}{c}
-\frac{\sqrt{3}}{2}\\\\\frac{1}{2}\end{array}\right)\quad
d_3=\left(\begin{array}{c}
\frac{\sqrt{3}}{2}\\\\\frac{1}{2}\end{array}\right)
\end{eqnarray*}
With these assumptions, one obtains the determinant of $M$:
\begin{eqnarray*}
|M|=6\sqrt{6}s^{11/4}+24\sqrt{3}s^3
\end{eqnarray*}
Similarly the determinant of $M$ for arbitrary angles is
\begin{eqnarray*}
|M|=32(\sin\theta_{13}\sin^2\theta_{12}+\sin\theta_{12}\sin^2\theta_{13})s^3-16\sqrt{2}(\sin\theta_{12}\sin\theta_{13}
\sin(\theta_{12}+\theta_{13}))s^{11/4}
\end{eqnarray*}
where $\theta_{12}$, $\theta_{13}$ are the angles between the
curves as shown in Fig.\ref{fig:fig2}. $M$ is nonsingular for any
$s$ with $\textrm{Re}(s)>0$ if $0<\theta_{12},\theta_{13} <\pi$.
The constraint on $\theta$ is not an issue since it has included
all the cases of interest.

\subsection{Analysis of Artificial Tangential Conditions}
As have been mentioned before, there are several options for the
artificial tangential conditions. We are interested to know if
different choices will lead to the same solution which is shown to
be true. To prove this point, it suffices to prove that the
position of the junction and the three tangential directions do
not depend on the artificial tangential conditions. The idea to
prove this point is to show that they all lead to the same
solution for $X_1$. If this is true, the position of the junction
point and the tangential direction of $X_1$ are uniquely
determined. Since the angle conditions are guaranteed, the
tangential directions of the other two curves could also be
uniquely determined. To sum up, the key point is proving
coefficients of $X^1$, i.e., $A_{11},A_{12}$ do not depend on the
extra conditions.

The coefficients $A_{ij},B_{ij}$ in solution (\ref{eq:sol}) can be
solved by
\begin{equation}
M \cdot C =P
\end{equation}
where $M$ is the coefficient matrix (\ref{eq:matrix}) for boundary
conditions, $C=[A_{11},A_{12},\cdots, A_{32},B_{32}]$ is the
coefficient vector to be solved and
$P=[p_1,p_2,\cdots,p_9,p_{10}]$ is a constant vector depending on
the initial data. Note that only $p_9,p_{10}$ and the last two
lines of $M$ depend on artificial tangential conditions.

According to the discussion above we need to prove $A_{11},A_{12}$
do not depend on the artificial tangential conditions. More
precisely, we need to prove $A_{11},A_{12}$ do not depend on the
last two lines of matrix $M$ and $p_9,p_{10}$.

For convenience, we rewrite $M$ into a block form
\begin{equation}
M=\left(\begin{array}{cc}
        M_1(8\times 2)& M_2(8\times 8)\\\\
        M_3(2\times 2) & M_4(2\times 8)
        \end{array}
        \right)
        \end{equation}
We do the Gauss elimination for block $M_2$ and it shows that the
rank of submatrix $M_2$ is 6 for any angle conditions. this means
we can make the last two lines of $M_2$ be zeros by row deduction
and meanwhile making the last two lines of $M_1$ into a full rank
$(2\times 2)$ matrix.

We again use $M$ to denote the new matrix after row deduction.
Next we compute $M^{-1}$ in a block form satisfying
\begin{eqnarray}\label{eq:times}
M\times M^{-1}&=&\left(\begin{array}{cc}
        M_1(8\times 2)& M_2(8\times 8)\\\\
        M_3(2\times 2) & M_4(2\times 8)
        \end{array}
        \right)\times
\left(\begin{array}{cc}
        \bar{M}_1(2\times 8)& \bar{M}_2(2\times 2)\\\\
        \bar{M}_3(8\times 8) & \bar{M}_4(8\times 2)
        \end{array}
        \right)\nonumber\\ \nonumber\\
        &=&\left(\begin{array}{cc}
        I(8\times 8)& \mathbf{0}\\\\
        \mathbf{0} & I(2\times 2)
        \end{array}
        \right)
        \end{eqnarray}
Expand directly to get
\begin{eqnarray}
&&M_1\times \bar{M}_1+M_2\times \bar{M}_3=\mathbf{I(8\times 8)} \label{eq:e1}\\
&&M_1\times \bar{M}_2+M_2\times \bar{M}_4=\mathbf{0(8\times
2)}\label{eq:e2}
\end{eqnarray}
Note that (\ref{eq:e1})-(\ref{eq:e2}) do not involve $M_3,M_4$
which means they do not depend on the artificial tangential
conditions. If $\bar{M}_1,\bar{M}_2$ can be determined by equation
(\ref{eq:e1})-(\ref{eq:e2}) then we can say $\bar{M}_1,\bar{M}_2$
do not depend on the artificial conditions. The fact
\begin{eqnarray*}
\left(\begin{array}{c}
        A_{11}\\
        A_{12}
        \end{array}
        \right)
        =\left(\begin{array}{cc}
        \bar{M}_1& \bar{M}_2
        \end{array}
        \right)\times
        P
        \end{eqnarray*}
implies that $A_{11},A_{12}$ do not depend on the artificial
conditions if we can further prove $\bar{M}_2=\mathbf{0}$ .

Actually, $\bar{M}_1$ can surely be solved from equation
(\ref{eq:e1}). This is because the last two lines of $M_2$ are
zeros and we have exactly sixteen equations involving only the
sixteen unknowns of $\bar{M}_1$. For the same reason we can solve
for $\bar{M}_2$ by equation (\ref{eq:e2}). Actually, since the
last two lines of $M_1$ is a full rank $(2\times 2)$ matrix
$\bar{M}_2$ must be zero. This completes the proof that the
coefficients $A_{11},A_{12}$ in equation (\ref{eq:linearized}) do
not depend on the artificial conditions. And consequently, the
shapes of the three curves do not depend on the artificial
tangential conditions. Novick-Cohen et al.~\cite{amy0} also
pointed out that the artificial conditions do not influence the
solutions, although the problem there is a little bit different.
In~\cite{amy0} the authors look at a three phase problem in which
all three interfaces evolve by minus the surface Laplacian of mean
curvature and meet at a triple junction.

\section{Numerical Discretization}\label{se:5}

Back to the full nonlinear problem, we present in detail the
discretization procedure of the parabolic scheme (\ref{eq:sys})
and junction conditions (\ref{eq:cond1})-(\ref{eq:mass}). The
basic approach is to use a staggered grid in $\sigma$ and we shall
denote the approximations by capital letters with subscripts,
i.e., $X_j(t)\simeq X((j-1/2)h,t)=(u((j-1/2)h,t),v((j-1/2)h,t))$
where $h$ is grid spacing and $N=1/h$ is the number of interior
grid points for $\sigma\in [0,1]$.

In order to write the discretized equations we introduce some
additional notations. Let $D_k$ denote the second order centered
approximation of the $k$th derivative, i.e.,
\begin{eqnarray*}
D_1X_j&=&(X_{j+1}-X_{j-1})/2h\\
D_2X_j&=&(X_{j+1}+X_{j-1}-2X_j)/h^2
\end{eqnarray*}
and let $D_+$ and $\mathcal{F}$ denote forward differencing and
forward averaging, respectively,
\begin{eqnarray*}
D_+X_j&=&(X_{j+1}-X_{j})/h\\
\mathcal{F}X_j&=&(X_{j+1}+X_{j})/2
\end{eqnarray*}
We discretize each motion separately.
\subsection{Grain Boundary Motion(Motion by Mean Curvature)}
The grain boundary motion is approximated at all grid points by
standard differences,
\begin{equation}
\dot{X_j^i}=\frac{D_2X_j^i}{|D_1X_j^i|^2} \qquad i=1,
j=1,2,\cdots,N
\end{equation}
where $\dot{X_j}$ stands for time derivative. Formally, these
discrete equations require values of $X_0$ and $X_{N+1}$ outside
the computation domain. We shall use the boundary condition to
extrapolate the interior values of $X_1$ and $X_N$ to the unknown
exterior values of $X_0$ and $X_{N+1}$. We shall give the details
of the extrapolation procedure later.

\subsection{Surface Diffusion}
The higher order derivatives appeared in surface diffusion are
approximated by
\begin{equation}
(X_{\sigma\sigma\sigma})_j\simeq
D_3X_j=\frac{D_2X_{j+1}-D_2X_{j-1}}{2h}
\end{equation}
\begin{equation}
(X_{\sigma\sigma\sigma\sigma})_j\simeq
D_4X_j=\frac{D_2X_{j-1}+D_2X_{j+1}-2D_2X_j}{h^2}
\end{equation}
There are some other terms such as $S_\sigma,S_{\sigma\sigma},
S_{\sigma\sigma\sigma}$ to be approximated. Start from
(\ref{eq:9}) and differentiate several times with respect to
$\sigma$ to get
\begin{eqnarray*}
S_\sigma&=&\sqrt{u_{\sigma}^2+v_{\sigma}^2}=|X_\sigma|\\
S_{\sigma\sigma}&=&\frac{u_{\sigma}u_{\sigma\sigma}+v_{\sigma}v_{\sigma\sigma}}{\sqrt{u_{\sigma}^2+v_{\sigma}^2}}
=\frac{X_\sigma\cdot X_{\sigma\sigma}}{|X_\sigma|}\\
S_{\sigma\sigma\sigma}&=&-\frac{(X_\sigma \cdot
X_{\sigma\sigma})^2}{|X_\sigma|^3}+ \frac{X_{\sigma\sigma}\cdot
X_{\sigma\sigma}+X_\sigma\cdot X_{\sigma\sigma\sigma}}{|X_\sigma|}
\end{eqnarray*}
Every term in scheme (\ref{eq:sys}) is now ready to be
approximated by standard differences.
\subsection{Junction Conditions at $\sigma=0$}
\begin{figure}
\begin{center}
\includegraphics[width = 10cm, height = 5.1cm, clip]{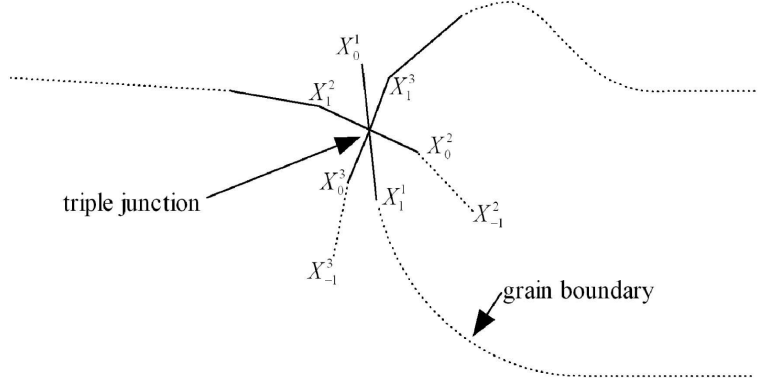}
\caption{Sketch of the ghost points at the triple junction.}
\label{fig:junc}
\end{center}
\end{figure}

The discretization at the junction point is much more complicated.
Since there are fourth order derivatives for the surface diffusion
we shall need two ghost points for each surface curve and one
ghost point for grain boundary. These ghost points are denoted by
$X^1_{0},X^2_{-1},X^2_{0},X^3_{-1}, X^3_{0}$ respectively. The
junction conditions (\ref{eq:cond1})-(\ref{eq:mass}) are
approximated as follows, see Fig. \ref{fig:junc}.

Condition (\ref{eq:cond1}),
\begin{equation}\label{eq:cond1d}
\mathcal{F}X^1_0=\mathcal{F}X^2_0=\mathcal{F}X^3_0=C
\end{equation}
where $C$ denotes the junction point.

The angle conditions (\ref{eq:cond2})-(\ref{eq:cond3}) are
approximated by
\begin{eqnarray}
\frac{D_+X_0^1}{|D_+X_0^1|}\cdot\frac{D_+X_0^2}{|D_+X_0^2|}&=&\cos\theta_{12}\\
\frac{D_+X_0^1}{|D_+X_0^1|}\cdot\frac{D_+X_0^3}{|D_+X_0^3|}&=&\cos\theta_{13}
\end{eqnarray}
Discretize condition (\ref{eq:potentials}) for each surface curve
to get
\begin{equation}
\frac{D_2X_C^2\cdot
(D_1X_C^2)^\perp}{|D_1X_C^2|^3}=-\frac{D_2X_C^3\cdot
(D_1X_C^3)^\perp}{|D_1X_C^3|^3}
\end{equation}
Since staggered grid are used, center $C$ is a midpoint not a grid
points. But we still can use previous notations $D_k$ with the
following extensions
\begin{eqnarray*}
X_{C-1}^i&=&(X_{-1}^i+X_0^i)/2=\mathcal{F}X_{-1}^i\\
X_{C+1}^i&=&(X_{1}^i+X_{2}^i)/2=\mathcal{F}X_{1}^i
\end{eqnarray*}
$\kappa_s$ can be expressed by
\begin{equation}
\kappa_s=\frac{X_{\sigma\sigma\sigma}\cdot
X_\sigma^\perp}{|X_\sigma|^4}-
3\frac{S_{\sigma\sigma}(X_{\sigma\sigma}\cdot X_\sigma^\perp)}
{|X_\sigma|^5}
\end{equation}
Thus condition (\ref{eq:mass}) is approximated by
\begin{equation}
\frac{D_3X^2_C\cdot (D_1X^2_C)^\perp}{|D_1X^2_C|^4}-
3\frac{S^2_{\sigma\sigma}(D_2X^2_C\cdot (D_1X^2_C)^\perp)}
{|D_1X^2_C|^5}=\frac{D_3X^3_C\cdot
(D_1X^3_C)^\perp}{|D_1X^3_C|^4}-
3\frac{S^3_{\sigma\sigma}(D_2X^3_C\cdot (D_1X^3_C)^\perp)}
{|D_1X^3_C|^5}
\end{equation}
Finally, the artificial tangential conditions is calculated by
\begin{equation}
\frac{D_2X_C^i\cdot D_1X_C^i}{|D_1X_C^i|^3}=0 \,\,\quad for
\,i=2,3
\end{equation}
We now finish discretizing the junction conditions.
\subsection{Domain Boundary Conditions at $\sigma=1$}
The discretization at $\sigma=1$ is straightforward.
\begin{eqnarray*}
&&\mathcal{F}X_N^i\big|_{t=n\cdot
dt}=\mathcal{F}X_N^i\big|_{t=(n-1)\cdot dt}\quad\textrm{for}\,
i=1,2,3 \\
&&D_+X_N^i=0 \quad \textrm{for} \,i=2,3
\end{eqnarray*}
\section{Time Stepping}\label{se:6}
\subsection{Explicit Scheme}
As an explicit scheme, forward Euler method is used for the time
stepping process.
\begin{equation}
X^{n+1}=X^n+\Delta tF(X^n)
\end{equation}
Here $F(X^n)$ denotes the right hand side in formulation
(\ref{eq:scheme}) evaluated at time level $n$. Time steps $\Delta
t$ are chosen so that the full discrete scheme is stable. Here we
choose $\Delta t=1e-12$. This scheme is easy to implement. Given
the results at time $n$ we update the values of the interior grid
points by forward Euler method to time level $n+1$ for the three
curves respectively. Next solve for the ghost points, junction
point and the boundary points by the boundary conditions. Then go
on to the next time level. The time step is excessively small due
to the stiffness of the fourth order parabolicity as noted
previously.
\subsection{Implicit Scheme}
In order to avoid the excessively small time steps for explicit
scheme we consider implicit techniques in this section. For
simplicity we use backward Euler method. Given the values at time
$n$ we update the values at time level $n+1$
 by solving the nonlinear system
\begin{equation}
X^{n+1}-\Delta tF(X^{n+1})-X^n=0
\end{equation}

Since the three curves are strongly coupled by the junction, we
solve all unknown points simultaneously including the ghost points
and the extrapolated boundary points. This leads to a large
nonlinear system which is solved by Newton's method. There is no
doubt that this scheme should be stable for any time steps. But it
can not survive a long time computation due to the nonuniform
tangential velocity which leads to a nonuniform distribution of
the grid points. This phenomenon can not be fixed even if we
refine the grid.

One way to overcome this difficulty is regridding the grid points
once they become too far or too close. But the bad distribution
could happen only near the junction and the closer to the junction
the sparser (or denser) the grid points are. Hence it is hard to
regrid no matter globally or locally. Another way is adjusting the
tangential velocity of the grid points such that they could adjust
themselves being uniform. And this is the motivation for the next
section.

\begin{figure}[h]
\setlength{\unitlength}{1cm}
\begin{picture}(10,6.5)
\put(0.5,0.5){\includegraphics[width=7cm,height=6cm]{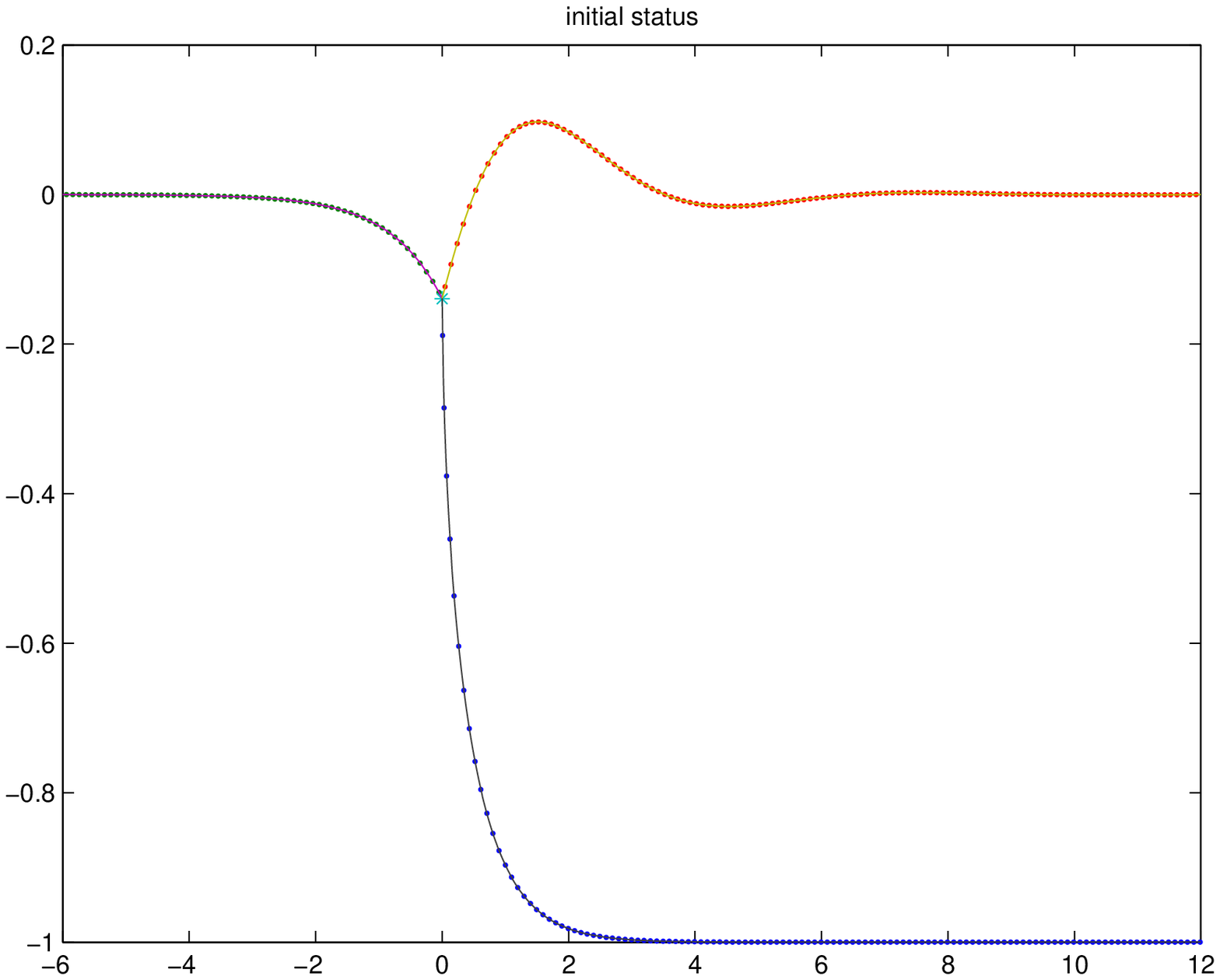}}
\put(8.5,0.5){\includegraphics[width=7cm,height=6cm]{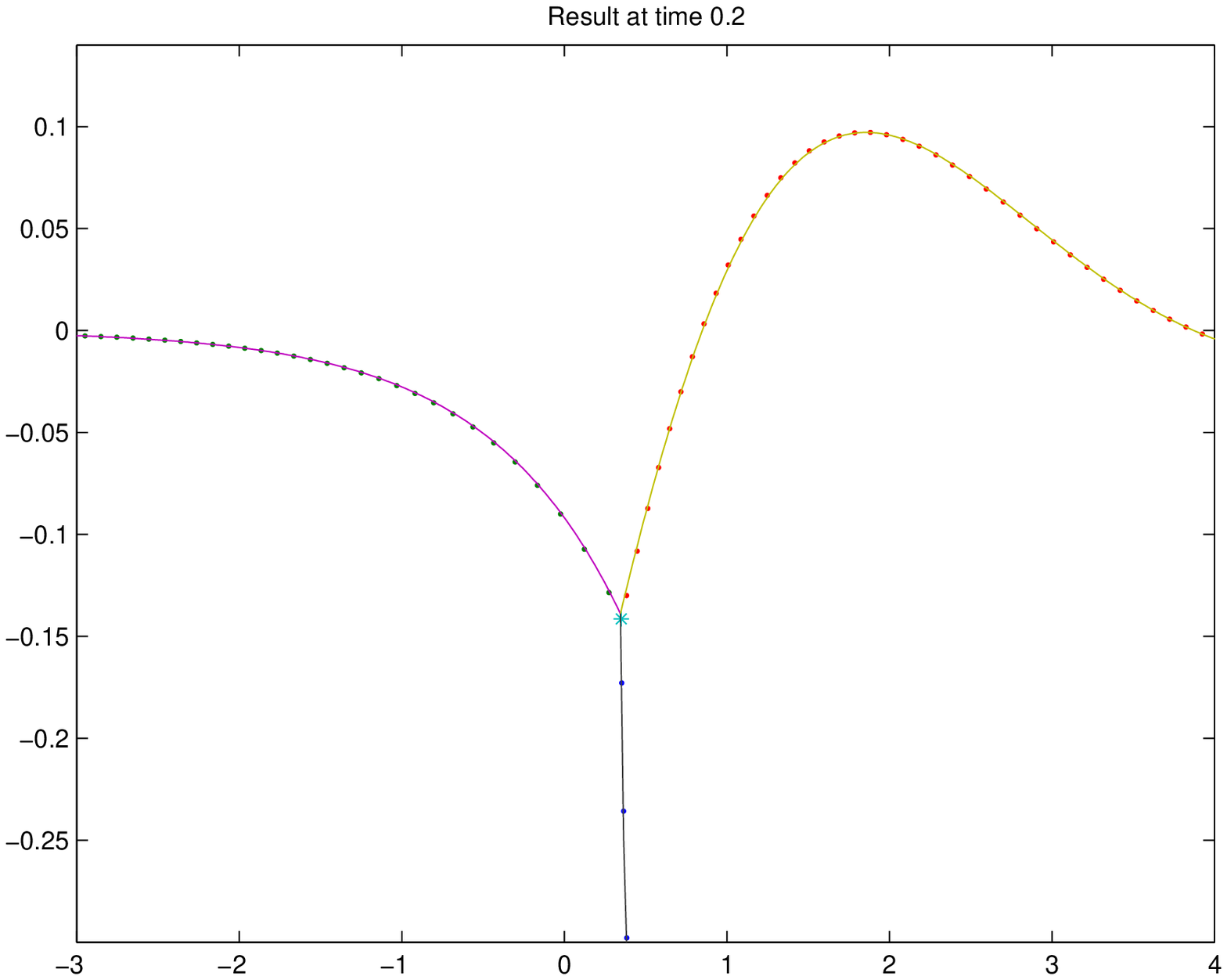}}
\end{picture}
\caption{Plot of results for scheme (\ref{eq:sys}) with backward
Euler method for a short time with $m=0.5$. Left: initial status
with grid points; Right: result zoomed in near triple junction.
Dotted line: numerical solution; Solid line: travelling wave
solution; Time step size: $\Delta t=1e-2$.}\label{fig:fig3}
\end{figure}

Numerical results for scheme (\ref{eq:sys}) with time step $\Delta
t=1e-4$ are shown in Fig. \ref{fig:fig3}. All numerical
experiments in this paper start from the same position as showed
in Fig. \ref{fig:fig3}. All results are compared with a travelling
wave solution solved by Amy et al~\cite{amy2}.
\\
\section{Adjustment of Tangential Velocity}\label{se:tangential}
We have mentioned that long time computations are problematical
even for implicit schemes. This is because of the bad distribution
of grid points. To get a more uniform distribution of gird points
along the curve we consider adding an artificial term to adjust
the tangential velocity of the grid points for the motion by
surface diffusion.

\begin{figure}[h]
\setlength{\unitlength}{1cm}
\begin{picture}(10,6.5)
\put(0.5,0.5){\includegraphics[width=7cm,height=6cm]{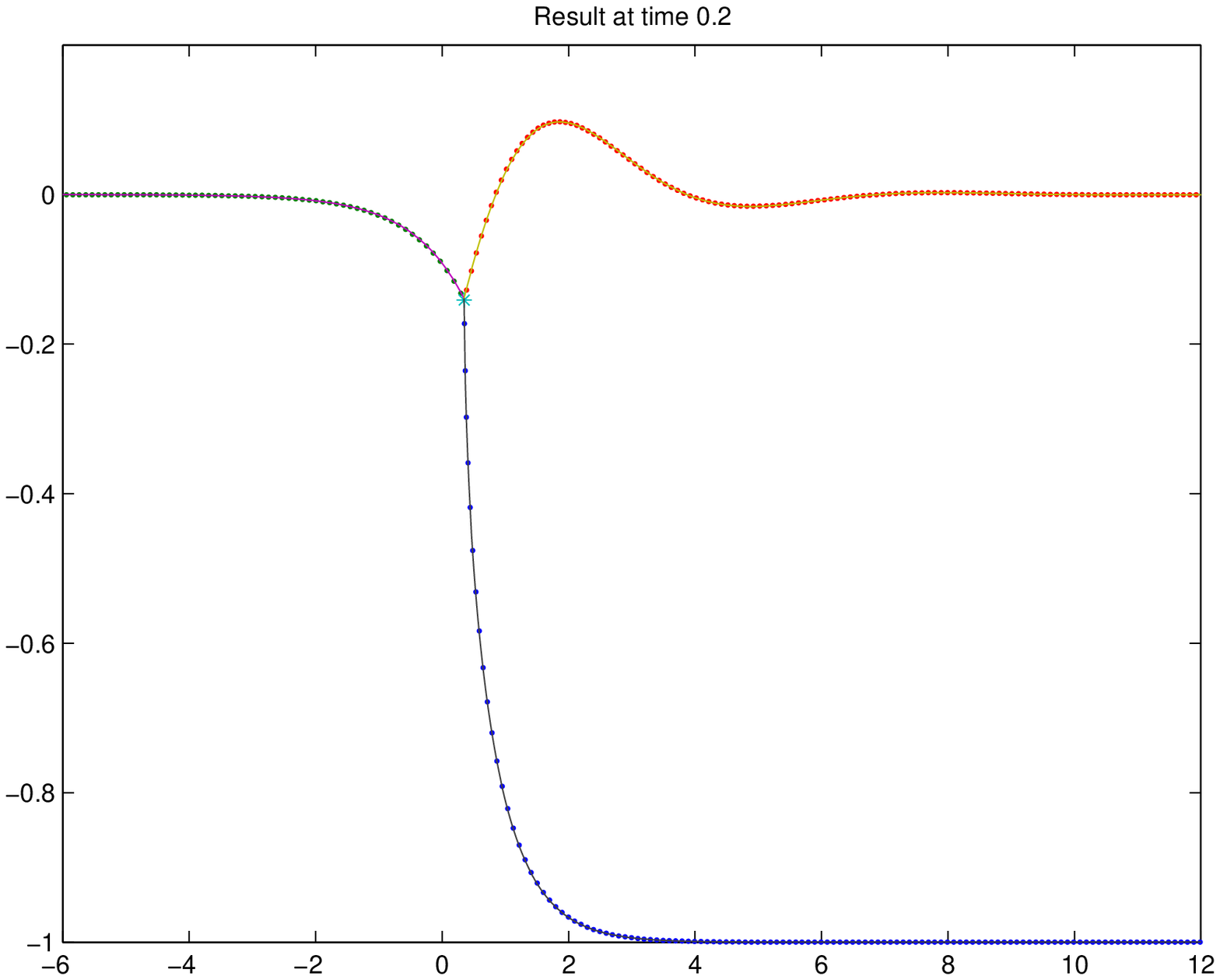}}
\put(8.5,0.5){\includegraphics[width=7cm,height=6cm]{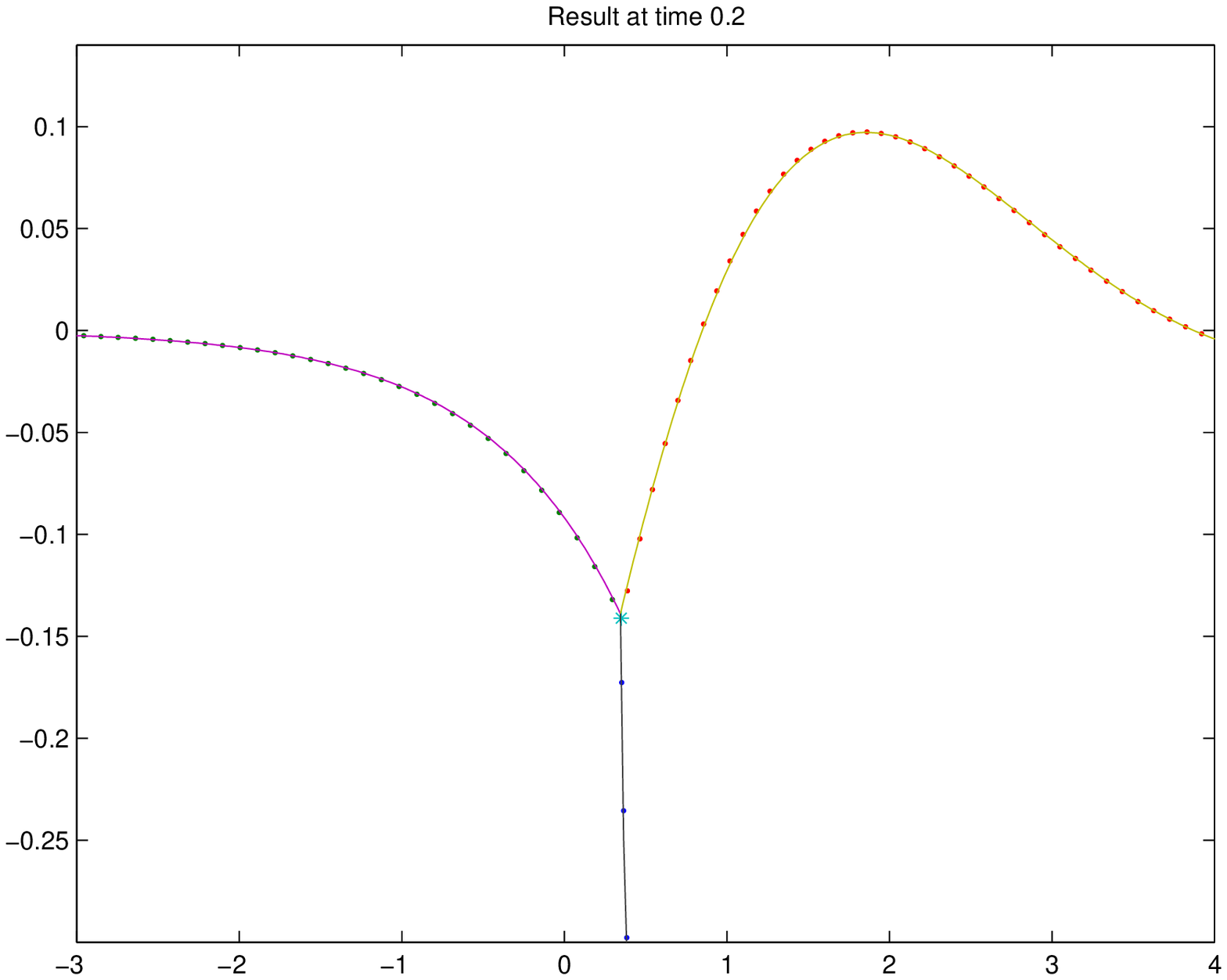}}
\end{picture}
\caption{Plot of the results for scheme (\ref{eq:newsc}) with
$m=0.5, \alpha=-100, \Delta t=0.01$. Left: result at t=0.2; Right:
result zoomed in near triple junction at t=0.2. Dotted line:
numerical result; Solid line: travelling wave solution.
}\label{fig:fig4}
\end{figure}

We consider the following modified scheme for the fourth order
problem
\begin{equation}\label{eq:newsc}
X_t=F(X)+\alpha
(\frac{X_{\sigma\sigma\sigma\sigma}}{|X_\sigma|^4}\cdot\vec{t})\,\vec{t}
\end{equation}
where $\alpha$ is a constant to be determined. The newly added
term in (\ref{eq:newsc}) will not influence the normal velocity
but it does change the tangential velocity. We do not know exactly
how to choose the optimal $\alpha$ but $\alpha=-100$ seems to work
well for our problem. The result is shown in Fig.\ref{fig:fig4}.
It is obvious that the grid points are much more uniform than that
in Fig.\ref{fig:fig3}. The time step size for Fig.\ref{fig:fig4}
is $\Delta t=0.01$. A numerical convergence study is shown in
Table \ref{tb:para}.
\begin{table}[!h]
\begin{center}
\begin{tabular}{|c||c||c|c|c|c|}
\hline
$dt$& $\Delta s$ & $L_2$ Norm & Rate & $L_\infty$ Norm& Rate\\
\hline\hline& 0.2 & 3.1494e-04 &&   2.0241e-03&\\
\cline{2-6}
 $dt=0.01\Delta s^2$  & 0.1 & 7.9775e-05 &1.9811 & 5.4797e-04&1.8852\\
\cline{2-6}
& 0.05 &2.1445e-05 &1.8953&1.4530e-05 &1.9151\\
\hline
\end{tabular}
\end{center}
\caption{Estimated errors and convergence rates for parabolic
formulation with $m=0.5$. Errors are evaluated at
$t=0.02$}.\label{tb:para}
\end{table}

Although the newly added tangential term improves the numerical
behavior it can not completely overcome the difficulty. The
artificial tangential conditions discussed in section \ref{se:4}
make the problem even more complicated. All these motivate us to
seek a more efficient scheme.
\\

\noindent\textbf{Remark} There is another way to adjust the
tangential velocity,
\begin{equation}
X_t=F(X)+\alpha
(\frac{X_{\sigma\sigma}}{|X_\sigma|^2}\cdot\vec{t})\,\vec{t}
\end{equation}
For this case, we should choose $\alpha$ positive, for example
$\alpha=100$.

\section{A PDAE Formulation}\label{se:8}
As we have pointed out in section \ref{se:tangential}, the fully
parabolic scheme (\ref{eq:sys}) does not always have good
numerical behavior. And the presence of the artificial tangential
condition makes the discretization of the original problem more
complicated. In this section we propose another formulation that
can overcome these disadvantages and also avoids possible loss of
tangential monotonicity in the parametrization due to the fourth
order PDE.

Let us again start from the motion by mean curvature. First of all
the basic evolution law should be satisfied, i.e.,
\begin{equation}\label{eq:curvature}
X_t\cdot \vec{n}-\kappa =0
\end{equation}
Because there are two free variables in this equation we need one
more equation for solvability. Since the requirement for the
normal direction motion has been fulfilled by equation
(\ref{eq:curvature}) we use the second equation to impose a
constraint on the distribution of grid points. It is natural to
let all grid points have equal spaces. To avoid introducing an
extra variable we let the change rate between any two adjacent
spaces is zero, i.e.,
\begin{equation}\label{eq:constraint}
|X_\sigma|_\sigma =0
\end{equation}
Note that
\begin{eqnarray*}
|X_\sigma|_\sigma =(\sqrt{X_\sigma\cdot
X_\sigma})_\sigma=\frac{X_\sigma\cdot
X_{\sigma\sigma}}{|X_\sigma|}
\end{eqnarray*}
the following equations are actually used to describe the motion
and keep grids equi-spaced,
\begin{eqnarray*}
&&X_t\cdot \vec{n}-\kappa =0\\
&&X_\sigma\cdot X_{\sigma\sigma}=0
\end{eqnarray*}
These are called partial differential algebraic equations (PDAEs).

In a similar way we derive the PDAEs for the motion by surface
diffusion,
\begin{eqnarray*}
&&X_t\cdot \vec{n}+\kappa_{ss} =0\\
&&X_\sigma\cdot X_{\sigma\sigma}=0
\end{eqnarray*}
Then the full PDAE system for the coupled motion is
\begin{eqnarray}\label{eq:daesys}
&&X^1_{t}\cdot \vec{n}-\kappa =0\nonumber\\
&&X^2_{t}\cdot \vec{n}+\kappa_{ss} =0\nonumber\\
&&X^3_{t}\cdot \vec{n}+\kappa_{ss} =0\\
&&X^i_{\sigma}\cdot X^i_{\sigma\sigma}=0 \qquad i=1,2,3\nonumber
\end{eqnarray}
The boundary conditions are the same as the parabolic case except
that we do not need artificial tangential conditions any more.

This is an implicit index-1 DAE system. Usually an index-1 DAE can
be discretized directly without any numerical difficulties
~\cite{ascher}, and that is our experience in this case.

Although the boundary conditions are the same as those of the
parabolic system, the discretization is a little bit different.
Instead of using five ghost points we now introduce only three
ghost points plus two extra variables which stand for the
curvature at the two ghost points corresponding to the two surface
branches. The two variables are denoted by $\kappa^2_0,\kappa^3_0$
and the last two junction conditions are approximated by
\begin{eqnarray*}
\frac{\kappa^2_0+\kappa^2_1}{2}&=&-\frac{\kappa^3_0+\kappa^3_1}{2}\\
\frac{(\kappa^2_0-\kappa^2_1)}{|D_1X^2_c|}&=&\frac{(\kappa^3_0-\kappa^3_1)}{|D_1X^3_c|}
\end{eqnarray*}
where $\kappa^i_1$ stands for the curvature of the first interior
point of curve $i$ and we use the average of $k^i_0, k^i_1$ to
approximate the curvature at the center point, i.e., junction
point. Again the sign should be carefully handled.

\begin{figure}[h]
\setlength{\unitlength}{1cm}
\begin{picture}(10,6.5)
\put(0.5,0.5){\includegraphics[width=7cm,height=6cm]{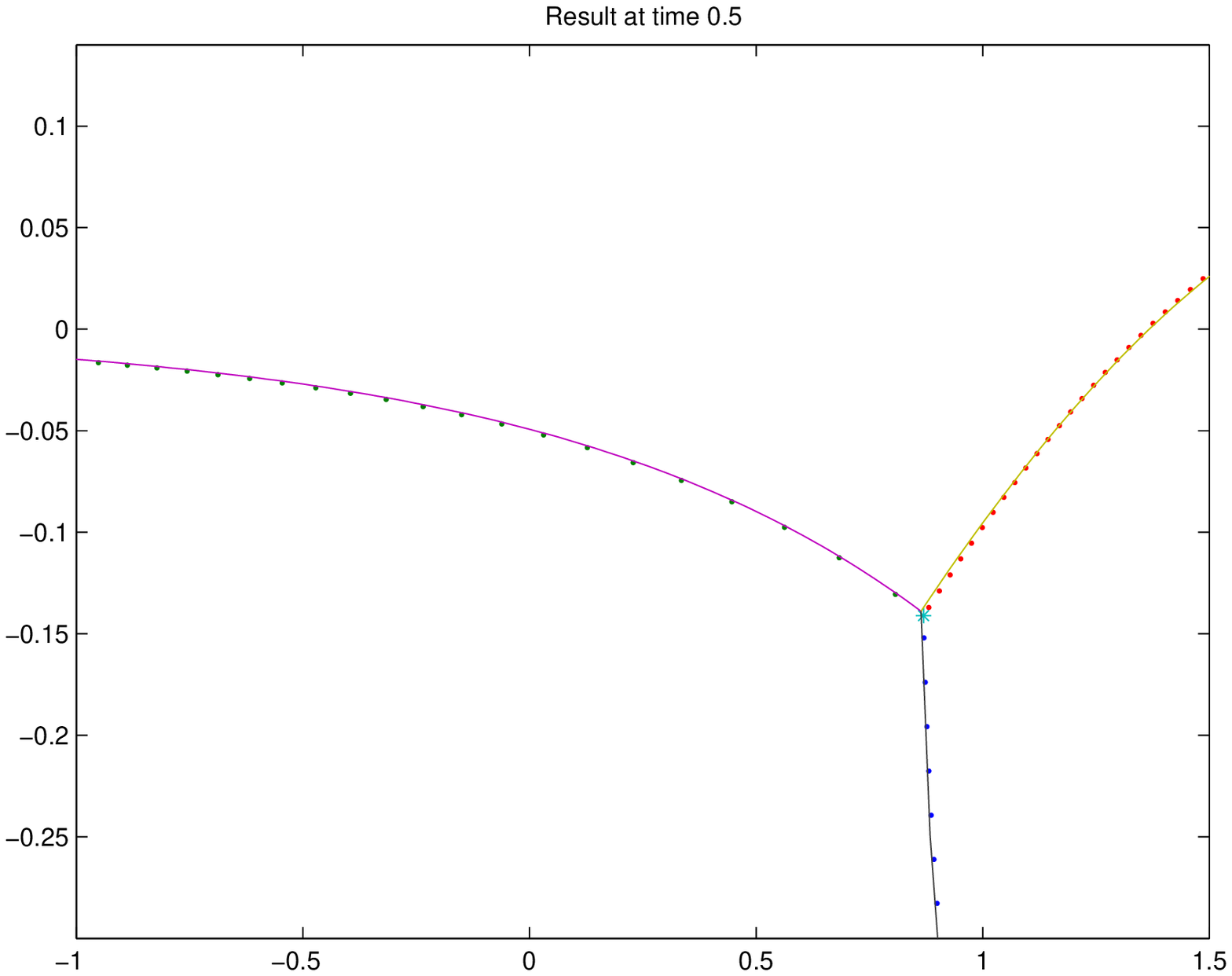}}
\put(8.5,0.5){\includegraphics[width=7cm,height=6cm]{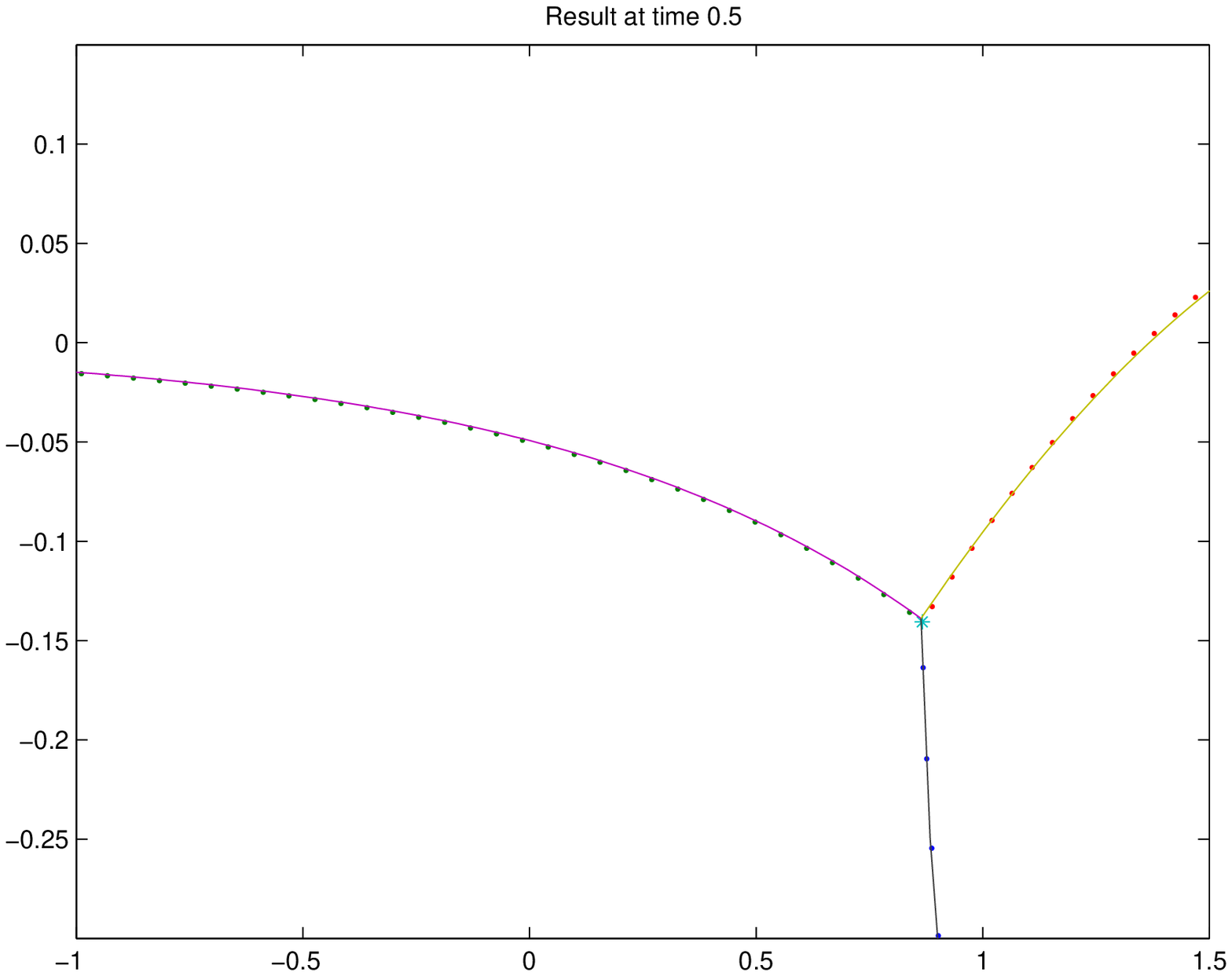}}
\end{picture}
\caption{Results comparison between the two schemes. Left: result
for (\ref{eq:sys}); Right: result for (\ref{eq:daesys}). Both
pictures are zoomed in near triple junction. Dotted line:
numerical solution; Solid line: travelling wave solution; Time
step size: $\Delta t=0.01$.}\label{fig:fig5}
\end{figure}

Implementing this scheme one obtains a better result shown in
Fig.\ref{fig:fig5}. The result is much more accurate and the grid
points are more uniform as well. An error comparison is shown in
Table \ref{tb:err}. A numerical convergence study of the PDAE
formulation is shown in Table \ref{tb:pdae}. The convergence rates
shown in Table \ref{tb:pdae} are close to 2 as expected.
\begin{table}[!h]
\begin{center}
\begin{tabular}{|c||c|c|}
\hline  &  Parabolic Formulation & PDAE Formulation\\
\hline\hline $\Delta s$ & 0.05 & 0.05 \\
\hline $\Delta t$ & 0.01 & 0.01 \\
\hline $L_\infty$ & 0.0041 & 0.0027
\\\hline
\end{tabular}
\end{center}
\caption{Performance of the two formulations with $L_\infty$ norm
and $\Delta t=0.01$.}\label{tb:err}
\end{table}
\begin{table}[!h]
\begin{center}
\begin{tabular}{|c||c||c|c|c|c|}
\hline
$dt$& $\Delta s$ & $L_2$ Norm & Rate & $L_\infty$ Norm& Rate\\
\hline\hline& 0.2 & 2.7837e-04 &&   1.8996e-03&\\
\cline{2-6}
 $dt=0.01\Delta s^2$  & 0.1 & 7.2717e-05 &1.9366 & 5.4444e-04&1.8029\\
\cline{2-6}
& 0.05 &1.8732e-05 &1.9568&1.4732e-04  &1.8858\\
\hline
\end{tabular}
\end{center}
\caption{Estimated errors and convergence rates for PDAE
formulation with $m=0.5$. Errors are evaluated at
$t=0.02$.}\label{tb:pdae}
\end{table}

\begin{figure}[h]
\setlength{\unitlength}{1cm}
\begin{picture}(10,6.5)
\put(4.5,0.5){\includegraphics[width=7cm,height=6cm]{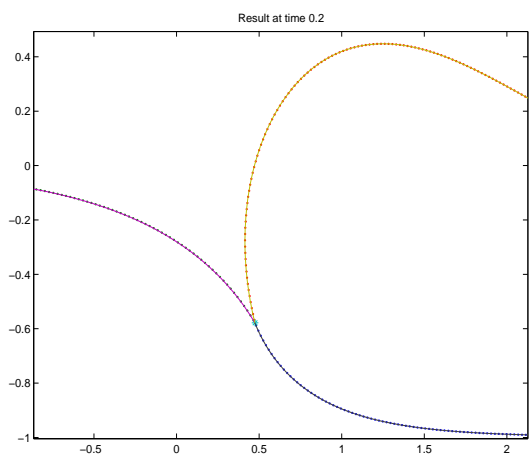}}
\end{picture}
\caption{Plot of the results for scheme (\ref{eq:daesys}) with
$m=1.96$ which has non-single valued upper surface. Dotted line:
numerical result; Solid line: travelling wave
solution.}\label{fig:nonsingle}
\end{figure}

Without difficulty we can apply this scheme to the case when the
surface curve is not a single-valued function as shown in
Fig.\ref{fig:nonsingle}
\\

\section{An Example of Surface Diffusion Problem}\label{se:star}
\begin{figure}[!h]
\begin{center}
\includegraphics[width = 11cm, height = 10cm, clip]{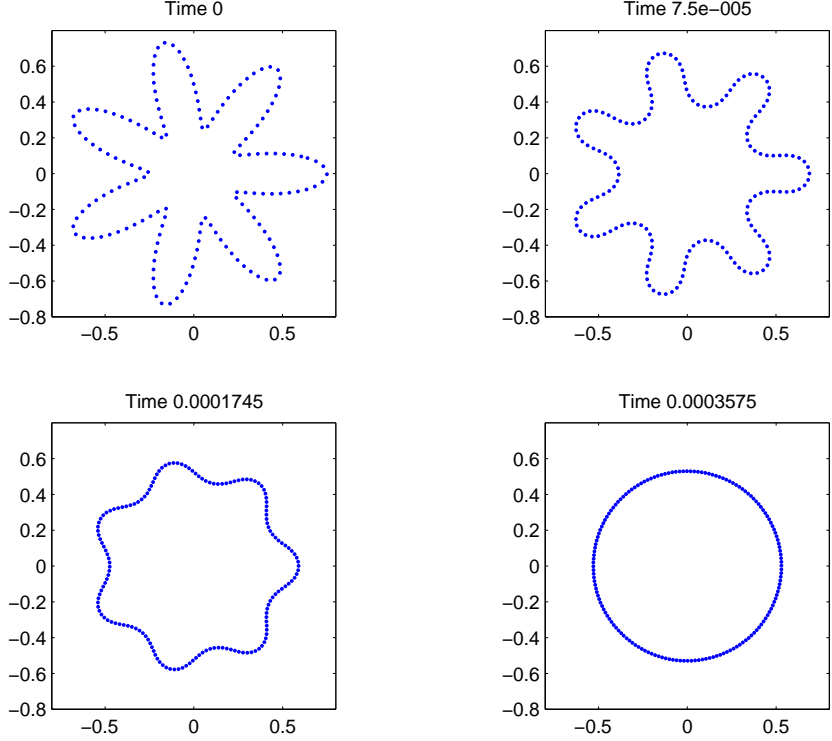}
\caption{A computational example that involve only the motion by
surface diffusion. $\Delta t=5\times 10^{-7}$. The area changes by
0.032\% } \label{fig:star}
\end{center}
\end{figure}

We temporally move our focus to a normal direction motion that
involves only motion by surface diffusion. The motion starts with
a closed star shaped curve and evolves with a speed equal to the
second derivative of curvature with respect to arc length.
According to the properties of surface diffusion the curve will
evolve into a circle and preserve the area enclosed by itself.
This problem is computed using the PDAE formulation for the
surface diffusion and the result is shown in Fig.\ref{fig:star}.
The method conserves the area quite well and the change is about
0.032\%. Similar examples have been investigated using level set
methods in~\cite{chopp,peter}. Level set methods have unbeatable
superiority for interface motion problem especially when there is
topology change. But for this simple problem (with no topology
change ) our method is more efficient and accurate. Note that this
problem can also be computed by the parabolic formulation.

\section{Well-posedness for the PDAE System}\label{se:9}
Similar to the parabolic system we do a well-posedness analysis
for the PDAE system in this section.

Considering the same linear problem as that in section
\ref{se:linearized} one obtains
\begin{eqnarray}
&&X^1_{t}\cdot d_1^\perp=X^1_{\sigma\sigma}\cdot d_1^\perp\nonumber\\
&& d_1\cdot X^1_{\sigma\sigma}=0\nonumber\\
&&X^2_{t}\cdot d_2^\perp=-X^2_{\sigma\sigma\sigma\sigma}\cdot d_2^\perp\nonumber\\
&& d_2\cdot X^2_{\sigma\sigma}=0\label{eq:2ndlinear}\\
&&X^3_{t}\cdot d_3^\perp=-X^3_{\sigma\sigma\sigma\sigma}\cdot d_3^\perp\nonumber\\
&& d_3\cdot X^3_{\sigma\sigma}=0\nonumber
\end{eqnarray}
where $d_i$ and $d_i^\perp$ stand for unit tangential direction
and unit normal direction of the $i^{th}$ curve respectively.

Linearization of the boundary conditions are exactly the same as
before. They differ only for the discretization procedure.

If $d_{i1},d_{i2}\neq0$ the linearized system (\ref{eq:2ndlinear})
has solution in the form {\large
\begin{eqnarray}\label{eq:1}
\left\{\begin{array}{l}
                    u_1=A_{11}e^{-\sqrt{s}\sigma}+B_{11}\\
                    v_1=-k_{1}A_{11}e^{-\sqrt{s}\sigma}+\frac{1}{k_{1}}B_{11}\\
                    u_2=A_{21}e^{\lambda_1\sigma}+B_{21}e^{\lambda_2\sigma}+C_{21}\\
                    v_2=-k_{2}(A_{22}e^{\lambda_1\sigma}+B_{22}e^{\lambda_2\sigma})+\frac{1}{k_{2}}C_{21}\\
                    u_3=A_{31}e^{\lambda_1\sigma}+B_{31}e^{\lambda_2\sigma}+C_{31}\\
                    v_3=-k_{3}(A_{32}e^{\lambda_1\sigma}+B_{32}e^{\lambda_2\sigma})+\frac{1}{k_{3}}C_{31}\\
                    \end{array}
                    \right.
\end{eqnarray}}
where $k_{i}=\frac{d_{i1}}{d_{i2}}$ is a constant and
\begin{eqnarray*}
\lambda_1=(-\frac{\sqrt{2}}{2}+\frac{\sqrt{2}}{2}i)\sqrt[4]{s}\quad
\lambda_2=(-\frac{\sqrt{2}}{2}-\frac{\sqrt{2}}{2}i)\sqrt[4]{s}
\end{eqnarray*}

Without changing the well-posedness of the problem we specify one
of the tangential directions, say $d_1=(0,-1)^T$. Further we
assume
\begin{eqnarray*}
\theta_{12},\theta_{13} \in (0,\pi) \textrm{ and
}\theta_{12},\theta_{13}\neq \frac{\pi}{2}
\end{eqnarray*}
Since $d_{11}=0$ now the solution is changed to
\begin{eqnarray}\label{eq:2}
\left\{\begin{array}{l}
                 u_1=A_{11}e^{-\sqrt{s}\sigma}\\
                 v_1=B_{11}\\
                  u_2=A_{21}e^{\lambda_1\sigma}+B_{21}e^{\lambda_2\sigma}+C_{21}\\
                    v_2=-k_{2}(A_{22}e^{\lambda_1\sigma}+B_{22}e^{\lambda_2\sigma})+\frac{1}{k_{2}}C_{21}\\
                    u_3=A_{31}e^{\lambda_1\sigma}+B_{31}e^{\lambda_2\sigma}+C_{31}\\
                    v_3=-k_{3}(A_{32}e^{\lambda_1\sigma}+B_{32}e^{\lambda_2\sigma})+\frac{1}{k_{3}}C_{31}\\
                    \end{array}
                    \right.
\end{eqnarray}
Apply these solutions to boundary conditions to get an $8 \times
8$ matrix $M$ and compute the determinant of $M$ directly to get
\begin{eqnarray*}
|M|=\frac{4\sqrt{2}s^{7/4}\sin(\theta_{12}+\theta_{13})-8s^2(\sin\theta_{12}+\sin\theta_{13})}
{\cos^2\theta_{12}\cos^2\theta_{13}}
\end{eqnarray*}
For the special case when one of the angles
$\theta_{12},\theta_{13}$ is $\frac{\pi}{2}$, for example,
$\theta_{12}=\frac{\pi}{2}$,
\begin{eqnarray}\label{eq:case2}
\left\{\begin{array}{l}
                 u_1=A_{11}e^{-\sqrt{s}\sigma}\\
                 v_1=B_{11}\\
                  u_2=C_{21}\\
                    v_2=A_{22}e^{\lambda_1\sigma}+B_{22}e^{\lambda_2\sigma}\\
                    u_3=A_{31}e^{\lambda_1\sigma}+B_{31}e^{\lambda_2\sigma}+C_{31}\\
                    v_3=-k_{3}(A_{32}e^{\lambda_1\sigma}+B_{32}e^{\lambda_2\sigma})+\frac{1}{k_{3}}C_{31}\\
                    \end{array}
                    \right.
\end{eqnarray}
The determinant of the coefficient matrix is
\begin{eqnarray*}
|M|=\frac{-8s^2\sin\theta_{13}+4\sqrt{2}s^{7/4}\cos\theta_{13}-8s^2}{\cos^2\theta_{13}}
\end{eqnarray*}
And if both $\theta_{12},\theta_{13}$ are $\frac{\pi}{2}$,
\begin{eqnarray}\label{eq:case2}
\left\{\begin{array}{l}
                 u_1=A_{11}e^{-\sqrt{s}\sigma}\\
                 v_1=B_{11}\\
                  u_2=C_{21}\\
                    v_2=A_{22}e^{\lambda_1\sigma}+B_{22}e^{\lambda_2\sigma}\\
                    u_3=C_{31}\\
                    v_3=A_{32}e^{\lambda_1\sigma}+B_{32}e^{\lambda_2\sigma}\\
                    \end{array}
                    \right.
\end{eqnarray}
The corresponding determinant of the matrix $M$ is
\begin{eqnarray*}
|M|=-16s^2
\end{eqnarray*}
Similar as before all cases discussed above are well-posed if
$0<\theta_{12},\theta_{13} <\pi$ and
$\theta_{12}+\theta_{13}\geq\pi$. Note that the well-posedness
property here coincides with that for the parabolic system we got
before.
\section{Conclusion}
We proposed two formulations to describe the coupled surface and
grain boundary motion. Both of them are well-posed and easy to be
implemented by finite difference method. Numerical results are
shown to be accurate. The PDAE formulation behaves better than the
parabolic form does. And since all grid points are equispaced for
the PDAE formulation it is convenient to regrid globally when
necessary. This often happens when the curves expand or shrink
quickly.

It is obvious that these schemes can also be used to simulate the
motion of a curve that involves only mean curvature motion or
surface diffusion as shown in section \ref{se:star}. And they are
extensible to any normal direction motion. Wherever applicable
these methods are more efficient comparing to the level set
methods. But they can not manage topology changes during the
evolution.

\section*{Acknowledgements}
The authors would like to thank Amy Novick-Cohen for the helpful
information on recent progress about the model addressed in this
paper.
\section*{Appendix A: Reformulation of the Motion.}
The original curvature motion and surface diffusion are given as
\begin{eqnarray*}
V_{c}&=&A\kappa\\
V_{d}&=&-B\kappa_{ss}
\end{eqnarray*}
Without loss of generality will shall prove in particular in a
parameterized form that we may normalize both $A$ and $B$ by
rescaling the time and space.

Suppose the original motions are modelled by
\begin{eqnarray}\label{eq:rescale}
X_t\cdot \vec{n}&=&A\kappa \nonumber\\
Y_t\cdot \vec{n}&=&-B\kappa_{ss}
\end{eqnarray}
We first rescale the time by
\begin{eqnarray*}
\tilde{t}=Tt
\end{eqnarray*}
Then equation (\ref{eq:rescale}) becomes
\begin{eqnarray}\label{eq:rescale2}
X_{\tilde{t}}\cdot \vec{n}&=&\frac{A}{T}\kappa \nonumber\\
Y_{\tilde{t}}\cdot \vec{n}&=&-\frac{B}{T}\kappa_{ss}
\end{eqnarray}
The spatial variable $X,Y$ are rescaled as
\begin{eqnarray*}
\tilde{X} = RX,\quad \tilde{Y} = RY
\end{eqnarray*}
And we can easily verify that
\begin{eqnarray*}
\tilde{s}=Rs
\end{eqnarray*}
where $\tilde{s}$ is the arc length in the new system. More
derivation shows that
\begin{eqnarray*}
\tilde{X}_{\tilde{s}}=X_{s}, \quad
\tilde{Y}_{\tilde{s}}=Y_{s},\quad
 \tilde{\kappa}=\frac{1}{R}\kappa, \quad \tilde{\kappa}_{\tilde{s}\tilde{s}}=\frac{1}{R^3}\kappa_{ss}\\
\end{eqnarray*}
Now equation (\ref{eq:rescale2}) becomes
\begin{eqnarray}\label{eq:rescale3}
\tilde{X}_{\tilde{t}}\cdot \vec{n}&=&\frac{AR^2}{T}\tilde{\kappa} \nonumber\\
\tilde{Y}_{\tilde{t}}\cdot
\vec{n}&=&-\frac{BR^4}{T}\tilde{\kappa}_{\tilde{s}\tilde{s}}
\end{eqnarray}
Here we use the same notation $\vec{n}$ since the normal direction
does not change. If we choose
\begin{eqnarray*}
R=\sqrt{\frac{A}{B}},\quad T=\frac{A^2}{B}
\end{eqnarray*}
then equation (\ref{eq:rescale3}) becomes
\begin{eqnarray}\label{eq:rescale4}
\tilde{X}_{\tilde{t}}\cdot \vec{n}&=&\tilde{\kappa} \nonumber\\
\tilde{Y}_{\tilde{t}}\cdot
\vec{n}&=&-\tilde{\kappa}_{\tilde{s}\tilde{s}}
\end{eqnarray}
This complete the normalization of coefficients $A$ and $B$.
\section*{Appendix B: Proof of Equation (\ref{eq:ksss}) }
One has the following fact
\begin{equation}\label{eq:kappa}
X_{ss}\cdot X_{ss}=\kappa^2
\end{equation}
Differentiating (\ref{eq:kappa}) with respect to $s$ one obtains
\begin{eqnarray*}\label{eq:26}
X_{sss}\cdot X_{ss}=\kappa\kappa_s
\end{eqnarray*}
\begin{eqnarray*}
X_{ssss}\cdot X_{ss}+X_{sss}\cdot X_{sss}
=\kappa\kappa_{ss}+\kappa_s^2
\end{eqnarray*}
Then an expression for $\kappa_{ss}$ can be derived,
\begin{equation}\label{eq:36}
\kappa_{ss}=\frac{X_{ssss}\cdot X_{ss}+X_{sss}\cdot X_{sss} -
\kappa_s^2}{\kappa}
\end{equation}
On the other hand, one has
\begin{eqnarray*}
X_{ss}=\kappa \vec{n}
\end{eqnarray*}
Take derivative again and calculate the inner product of $X_{sss}$
\begin{eqnarray*}
X_{sss}=\kappa_s\vec{n}+\kappa
\vec{n}_s=\kappa_s\vec{n}-\kappa^2\vec{t}
\end{eqnarray*}
\begin{eqnarray}\label{eq:39}
X_{sss}\cdot X_{sss} &=&\kappa_s^2+\kappa^4
\end{eqnarray}
Substitute equation (\ref{eq:39}) into (\ref{eq:36}) to get
\begin{equation}\label{eq:40}
\kappa_{ss}=X_{ssss}\cdot\vec{n}+\kappa^3
\end{equation}
Using equation (\ref{eq:40}) together with the fact that
\begin{equation}\label{eq:44}
\kappa^2X_{ss}\cdot \vec{n}=\kappa^2(X_{ss}\cdot \vec{n})=\kappa^3
\end{equation}
one obtains
\begin{eqnarray*}
\kappa_{ss}=(X_{ssss}+\kappa^2X_{ss})\cdot \vec{n}
\end{eqnarray*}
which completes the proof.


\end{document}